\documentclass{amsart}
\usepackage{graphicx}

\theoremstyle{plain}
\newtheorem{thm}{Theorem}[section]
\newtheorem{lemma}[thm]{Lemma}
\newtheorem{cor}[thm]{Corollary}
\newtheorem{propn}[thm]{Proposition}
\newtheorem{add}[thm]{Addendum}

\theoremstyle{definition}
\newtheorem{defn}[thm]{Definition}

\numberwithin{equation}{section}

\newcommand{\reals}{\ensuremath{\mathbb{R}}}
\newcommand{\half}{\ensuremath{\frac{1}{2}}}
\newcommand{\into}{\ensuremath{\hookrightarrow}}
\newcommand{\tb}{\mathop{\rm tb}\nolimits}

\newcommand{\id}{\mathop{\rm id}\nolimits}

\begin{document}

\title[Symplectic $2$-handles]{Symplectic $2$-handles and transverse links}

\author{David T. Gay}
\address{Department of Mathematics\\
University of Arizona\\
617 North Santa Rita\\
Post Office Box 210089\\
Tucson, Arizona 85721\\
U.S.A.}
\email{dtgay@math.arizona.edu}
\urladdr{http://www.math.arizona.edu/{\textasciitilde}dtgay}

\subjclass[2000]{Primary 57R17, 57R65; Secondary 57M99}
\keywords{Symplectic handle, symplectic handlebody, contact surgery,
fibered link, fibered knot, symplectic filling, convexity, concavity,
transverse link, transverse knot, symplectic germ}

\begin{abstract}

A standard convexity condition on the boundary of a symplectic
manifold involves an induced positive contact form (and contact
structure) on the boundary; the corresponding concavity condition
involves an induced negative contact form. We present two methods of
symplectically attaching $2$-handles to convex boundaries of
symplectic $4$-manifolds along links transverse to the induced contact
structures. One method results in concave boundaries and depends on a
fibration of the link complement over $S^1$; in this case the handles
can be attached with any framing {\em larger than} a lower bound
determined by the fibration. The other method results in a weaker
convexity condition on the new boundary (sufficient to imply tightness
of the new contact structure), and in this case the handles can be
attached with any framing {\em less than} a certain upper bound. These
methods supplement methods developed by Weinstein and Eliashberg for
attaching symplectic $2$-handles along Legendrian knots.

\end{abstract}

\maketitle

\section{Results and Motivation}
\label{S:ResMot}

When constructing symplectic manifolds it is natural to wonder whether
topological techniques using handles can be made to work
symplectically. Weinstein~\cite{Weinstein} and
Eliashberg~\cite{EliashStein} have shown how to do this in certain
cases; here we present two new symplectic ``handle-by-handle''
constructions in dimension four.

In such constructions it is desirable to retain control of the
symplectic form near the boundary; one form of control is the
following: Given a symplectic manifold $(X,\omega)$ we say that
$\partial X$ is {\em convex} (respectively {\em concave}) if there
exists a vector field $V$ defined in a neighborhood of $\partial X$,
satisfying the equation $\mathcal{L}_V \omega = \omega$ (in other
words, $V$ is a symplectic dilation) and pointing {\em out of}
(respectively {\em into}) $X$. This induces a contact form $\alpha =
\imath_V \omega|_{\partial X}$ and a contact structure $\xi = \ker
\alpha$ on $\partial X$.  Weinstein and Eliashberg show that, if
$(X,\omega)$ is a symplectic $2n$-manifold with $\partial X$ convex,
then one can attach $k$-handles to $X$, for $0 \leq k \leq n$, and
extend $\omega$ across the handles so that the new boundary is again
convex. Conditions are imposed on the attaching spheres in relation to
the contact structure $\xi$ on $\partial X$ and in particular, in
dimension four, $2$-handles must be attached along {\em Legendrian}
knots (knots tangent to $\xi$). In this paper we show how to
symplectically attach $2$-handles along {\em transverse} knots
(transverse to $\xi$) in the convex boundary of a symplectic
$4$-manifold so that the new boundary becomes {\em concave}. Along the
way, we see boundaries which are partially convex and partially
concave, so we develop a careful theory for such boundaries.

For a weaker form of control, we say that $\partial X$ is {\em
weakly convex} if $\partial X$ supports a positive contact structure
$\xi$ such that $\omega|_\xi$ is nondegenerate; convexity implies weak
convexity. In this paper we also show how to symplectically attach
$2$-handles along transverse knots in the convex boundary of a
symplectic $4$-manifold so that the new boundary becomes weakly convex.

\subsection{Theorems and Definitions}
We first present the statements of two theorems followed by the
definitions needed to understand the statements. Let $(X,\omega)$ be a
compact symplectic $4$-manifold, suppose that $\partial X$ is convex
with induced contact structure $\xi$ and let $L$ be a transverse link
in $(\partial X,\xi)$ with a chosen framing $F$.
\begin{thm} \label{T:WeakCvx}
If $L$ is {\em fat with respect to $F$} (a condition which holds for
negative enough framings), then there exists a symplectic $4$-manifold
$(Y,\eta)$ containing $(X,\omega)$, obtained by {\em enlarging}
$(X,\omega)$ in a neighborhood of $L$ and then attaching $2$-handles
along $L$ with framing $F$, such that $\partial Y$ is weakly convex.
\end{thm}
\begin{thm} \label{T:cvx2ccv}
If $L$ is {\em nicely fibered} and if $F$ is {\em positive with
respect to the fibration} (a condition which holds for positive enough
framings), then there exists a symplectic $4$-manifold $(Y,\eta)$
containing $(X,\omega)$, obtained by {\em enlarging} $(X,\omega)$ and
then attaching $2$-handles along $L$ with framing $F$, such that
$\partial Y$ is concave.
\end{thm}
(For both theorems, more precise statements will be made later
about the resulting contact structures on $\partial Y$.) We will now
define ``fat with respect to a framing'', ``nicely fibered'' and
``positive with respect to a fibration'' and we will say what we mean
by ``enlarging $(X,\omega)$''.

First we establish orientation conventions: In this paper all
manifolds are assumed to be oriented, and all symplectic forms are
assumed to be compatible with the given orientations. A vector field
$V$ in a manifold $X$ transverse to a codimension one submanifold $M$
is {\em positively} transverse if $V(p)$ followed by an oriented basis
for $T_p M$ is an oriented basis for $T_p X$. Boundaries are oriented
so that outward normals are positively transverse. Noncritical level
sets of a Morse function are oriented so that upward pointing vectors
are positively transverse. When $M$ is an oriented manifold, $-M$ will
refer to $M$ with the opposite orientation. Unless stated otherwise,
all diffeomorphisms and embeddings preserve orientations.

The first definitions involve local control in a neighborhood of a
transverse knot $K$.  For this it is convenient to use ``polar
coordinates'' near $K$, by which we mean functions $(r,\mu,\lambda)$
defining an embedding of a neighborhood $\nu$ of $K$ into $\reals^2
\times S^1$, with $(r,\mu)$ mapping to polar coordinates on $\reals^2$
and $\lambda$ mapping to $S^1$, and such that $K = \{r=0\}$. Note that
the function $\mu:\nu \setminus K \to S^1$ defines a framing of $K$,
which we label $F_\mu$. If $K$ is transverse to a contact structure
$\xi$, then we can always find {\em normal} polar coordinates on a
neighborhood $\nu$, by which we mean polar coordinates
$(r,\mu,\lambda)$ such that $\xi|_{\nu \setminus K}$ is spanned by
$\partial_r$ and $\pm \log(r^2) \partial_\mu +
\partial_\lambda$. (Where $\pm$ means $+$ when $\xi$ is positive and
$-$ when $\xi$ is negative.)  This follows from Darboux's theorem for
contact structures. (See, for example,~\cite{McDSal}.) A consequence
of this parametrization of $r$ is that, for each integer $n$, the
characteristic foliation on the torus $\{r=e^n\}$ is a family of
parallel longitudes realizing the framing $F_\mu + n$. Otherwise, the
exact parametrization of $r$ is not particularly important and so we
say that $(r,\mu,\lambda)$ is an {\em almost normal} coordinate system
if there exists a function $R$ such that $(R \circ r,\mu,\lambda)$ is
a normal coordinate system.  (If $(r,\mu,\lambda)$ is an almost normal
coordinate system then so is $(r,\mu - k \lambda, \lambda)$ for any $k
\in \mathbb{Z}$ and the relation between the framings is that
$F_{\mu-k\lambda} = F_\mu + k$.)

\begin{defn} \label{D:fat}
Given a knot $K$ in a contact $3$-manifold $(M,\xi)$ with a framing $F$
and with a chosen neighborhood $\nu$ with normal coordinates
$(r,\mu,\lambda)$, we say that the coordinate system {\em goes out as
far as} $F$ if the image of $\nu$ in $\reals^2 \times S^1$ contains
the solid torus $\{ r \leq e^n\}$, where $F = F_\mu + n$. Given a link
$L$ in $(M,\xi)$ and a framing $F$ of $L$, with $F_i$ being the
framing induced on each component $K_i$ of $L$, we say that $L$ is
{\em fat with respect to $F$} if there exist mutually disjoint
neighborhoods $\nu_i$ of each $K_i$, each with normal coordinates
which go out as far as $F_i$.
\end{defn}

A link $L$ in a $3$-manifold $M$ is fibered if $M \setminus L$ fibers
over $S^1$, and one often requires some controlled behavior of the
fibration near $L$. For each fiber $\Sigma$, a contact structure $\xi$
on $M$ induces a (singular) {\em characteristic foliation} on
$\Sigma$, the integral curves of $T\Sigma \cap \xi$. 

A {\em contact vector field} is a vector field which preserves a given
contact structure; we will say that a contact vector field is {\em
transverse} if it is everywhere transverse to the contact plane
field. We will be interested in transverse contact vector fields which
are also transverse to the fibers of a fibration.

Such a vector field $V$ co-orients both the fibers and the plane field
and hence orients the characteristic foliation on each fiber according
to the following convention: At a point in $M$ where $\xi$ is
transverse to a fiber, let $\alpha$ be a 1-form with kernel $\xi$ such
that $\alpha(V) > 0$ and let $\beta$ be a 1-form with kernel tangent
to the fiber such that $\beta(V) > 0$. Then let $\gamma$ be a 1-form
such that $\alpha \wedge \beta \wedge \gamma > 0$. Then $\ker \alpha
\cap \ker \beta$ is transverse to $\ker \gamma$ so we restrict
$\gamma$ to define the orientation on the characteristic
foliation. This orientation does not change if we replace $V$ with
$-V$.
\begin{defn} \label{D:nicefiber}
A transverse link $L = K_1 \cup \ldots \cup K_n \subset (M,\xi)$ is
{\em nicely fibered} if there exists a fibration $p: M \setminus L \to
S^1$ and a transverse contact vector field $V$ (defined on all of $M$)
satisfying the following conditions:
\begin{itemize}
\item
$V$ is transverse to the fibers of $p$.
\item
For each $K_i$ there exist almost normal polar coordinates
$(r,\mu,\lambda)$ on a neighborhood $\nu_i$ near $K_i$ such that, on
$\nu_i$, $\partial_r$ is tangent to the fibers, $dr(V)=0$ and $V$ and
$dp$ are both invariant under the flows of $\partial_{r}$,
$\partial_{\mu}$ and $\partial_{\lambda}$.
\item 
Letting $V$ co-orient both $\xi$ and the fibers, the oriented
characteristic foliation on the fibers near each $K_i$ should point in
towards $K_i$. (By the previous condition the unoriented foliation
will be by radial lines in each normal coordinate system.)
\end{itemize}
\end{defn}
Such a fibration has a bearing on framings of $L$. For each $K_i$ with
normal polar coordinates $(r,\mu,\lambda)$ as in the definition,
consider a torus $T=\{r=R\}$ with coordinates $(\mu,\lambda)$. The
fibers of $p$ intersect $T$ in a family of parallel lines with
``slope'' $\frac{d\mu}{d\lambda} = s_p \in \reals \cup \{\infty\}$. A
framing $F$ of $L$ also gives a family of parallel longitudes in $T$
with ``slope'' $\frac{d\mu}{d\lambda} = s_F \in \mathbb{Z}$. (In terms
of the earlier notation, $F=F_\mu+s_F$.)
\begin{defn} \label{D:posframing}
In this situation we say that $F$ is {\em positive with respect to the
fibration} if, for each $K_i$ as above, $s_p \neq \infty$ and $s_F > s_p$.
\end{defn}

\begin{defn} \label{D:enlarge}
By {\em enlarging} a symplectic manifold $(X,\omega)$ we will mean the
following process: Choose a smooth function $h:\partial X \to
[0,\infty)$ and let $C=\{(t,p) \mid 0 \leq t \leq h(p)\} \subset
\reals \times \partial X$. Then attach $C$ to $X$ by the obvious
identification of $\{0\} \times \partial X \subset C$ with $\partial X
\subset X$. Given a neighborhood $N$ in $\partial X$, we will say that
the enlargement is {\em supported inside $N$}, or simply that we are
{\em enlarging $(X,\omega)$ inside $N$}, if the support of $h$ is
contained in $N$. After attaching $C$ to $X$ choose some extension of
$\omega$ to a symplectic form on all of $X \cup C$ and replace
$(X,\omega)$ with the new, larger symplectic manifold. There is a
natural identification between the old $\partial X = \{0\} \times X$
and the new $\partial X = \{h(p),p\} \in C$.
\end{defn}
When $\partial X$ is convex there is a canonical extension of $\omega$
over $C$, constructed as follows: Given a contact form $\alpha$ on a
$3$-manifold $M$ there is a canonical symplectification $(\reals \times
M, d(e^t \alpha))$, in which $\partial_t$ is a symplectic dilation
inducing $\alpha$ on $\{0\} \times M$. If $M$ is embedded in another
symplectic manifold $(Y,\eta)$ with a symplectic dilation $V$
positively transverse to $M$ inducing $\alpha$, then flow along $V$
generates a symplectomorphism from a neighborhood of $\{0\} \times M$
in $(\reals \times M, d(e^t \alpha))$ to a neighborhood of $M$ in
$(Y,\eta)$, restricting to the identity on $M$ and sending
$\partial_t$ to $V$. Thus if $\alpha$ is the induced contact form on
$\partial X$ and we use the form $d(e^t \alpha)$ on $C \subset \reals
\times \partial X$, this will patch together smoothly with $\omega$ on
$X$ and the symplectic dilation will extend across $C$ to be
transverse to the new $\partial X$. Furthermore the induced form on
the new $\partial X$ is $e^h \alpha$, and so the underlying contact
structure is unchanged. For future reference, we will call this
enlargement the convex enlargement of height $h$; it is clearly useful
if we would like to rescale the contact form on $\partial X$ by some
function greater than or equal to $1$.

This is the enlargement used in the two theorems, although in
theorem~\ref{T:cvx2ccv} we will need to see the enlargement in a more
general setting. Enlarging $(X,\omega)$ does not change the
diffeomorphism type of $X$ and so can instead be thought of as a
deformation of $\omega$ keeping $X$ fixed.

\subsection{Discussion}
If $L$ is fat with respect to a particular framing then it is also fat
with respect to any more negative framing, so theorem~\ref{T:WeakCvx}
gives a construction that works for very negative framings but is less
likely to work the more positive the framings become.  If $F$ is
positive with respect to a nice fibration then any more positive
framing is also positive with respect to the fibration, so that
theorem~\ref{T:cvx2ccv} gives a construction that works for
very positive framings (assuming the fibration exists) but is less
likely to work the more negative the framings become.

It is interesting to see these constructions alongside the
construction of Weinstein~\cite{Weinstein} mentioned earlier. (We use
Weinstein as our source because Weinstein's discussion is strictly
symplectic whereas Eliashberg~\cite{EliashStein} discusses the
construction in the case of Stein manifolds.) To simplify matters we
present the result only in dimension four.
\begin{thm}[Weinstein] \label{T:weinstein}
Let $(X,\omega)$ be a symplectic $4$-manifold with $\partial X$ convex
with induced contact structure $\xi$. Then, given any Legendrian knot
$K \subset (\partial X,\xi)$ there exists a symplectic $4$-manifold
$(Y,\eta)$ containing $(X,\omega)$, with $\partial Y$ convex, obtained
by enlarging $(X,\omega)$ in a neighborhood of $K$ and then attaching
a $2$-handle along $K$ with framing $\tb(K)-1$ (where $\tb(K)$ is the
Thurston-Bennequin framing of $K$). We can also symplectically attach
any number of 1-handles to $(X,\omega)$ to get $(Y,\eta)$ with convex
boundary.
\end{thm}
(The Thurston-Bennequin framing of $K$ is the framing given by any
vector field transverse to $K$ but lying in $\xi$.)  Since any knot is
$C^0$-close to a Legendrian knot and every Legendrian knot can be
isotoped so as to make its Thurston-Bennequin framing more negative
(see~\cite{EliashKnots}), the $2$-handle part of this theorem is also a
construction which, given a smooth knot type, works for very negative
framings but is less likely to work the more positive the framings
become.

This result and the fact that the contact structure on a weakly convex
boundary is always tight~\cite{EliashFillHolDisk} were used by
Gompf~\cite{Gompf} to construct many $3$-manifolds with tight contact
structures, beginning with the standard positive contact structure on
$S^3$ as the convex boundary of $B^4$ with its standard symplectic
structure. Consider the following observation (suggested by John
Etnyre):
\begin{propn} \label{P:LegXvrs}  
  Suppose that $K$ is a Legendrian knot in a positive contact
  $3$-manifold with a given neighborhood $\nu$ and a framing $F \leq
  \tb(K) - 1$. Then there exists a transverse knot $K^\prime$ inside
  $\nu$, isotopic to $K$, which is fat (inside $\nu$) with respect to
  $F$.
\end{propn}
This tells us that in fact theorem~\ref{T:WeakCvx}, together with the
1-handle part of theorem~\ref{T:weinstein}, can be used to construct
the same $3$-manifolds that Gompf constructs, also with tight contact
structures, but now they are weakly convex rather than convex
boundaries of symplectic $4$-manifolds.

Constructing manifolds with concave boundaries is interesting for two
reasons. First we get some answers to a simple symplectic filling
question: which $3$-manifolds with which contact structures can be
realized as the concave boundaries of symplectic $4$-manifolds?
Secondly, if we can carefully characterize the contact structures that
result from our concave constructions, and also construct symplectic
$4$-manifolds with convex boundaries and carefully characterize the
resulting contact structures, then we may be able to use the following
standard glueing construction to produce closed symplectic
$4$-manifolds: Let $(X_1,\omega_1)$ and $(X_2,\omega_2)$ be symplectic
$4$-manifolds with $\partial X_1$ convex and $\partial X_2$ concave
with induced contact forms $\alpha_1$ and $\alpha_2$ and contact
structures $\xi_1$ and $\xi_2$. Suppose that $\partial X_1$ and
$\partial X_2$ are connected and that there exists a contactomorphism
$\phi: (\partial X_1,\xi_1) \to (-\partial X_2,\xi_2)$. Then $\phi^*
\alpha_2 = f \alpha_1$ for some nonzero function $f$ on $\partial
X_1$. Rescaling $\omega_2$ by a constant we can arrange that $f >
1$. Let $h = \log f$. Enlarging $(X_1,\omega_1)$ with the canonical
convex enlargement of height $h$ will then replace $\alpha_1$ with
$e^h \alpha_1 = f \alpha_1$. Thus we can in fact arrange that $\phi^*
\alpha_2 = \alpha_1$. Finally, $(X_1,\omega_1)$ can be glued to
$(X_2,\omega_2)$ by identifying a neighborhood of $\partial X_1$ with
$((-\epsilon,0] \times \partial X_1, d(e^t \alpha_1))$ and, using
$\phi$, identifying a neighborhood of $\partial X_2$ with
$([0,\epsilon) \times \partial X_1, d(e^t \alpha_1))$.

For specific examples where theorem~\ref{T:cvx2ccv} applies, observe
that the transverse unknot and the Hopf link in $S^3$ with the
standard contact structure are both nicely fibered, with framings
greater than $0$ being positive with respect to the unknot fibration
and framings greater than or equal to $0$ being positive with respect
to the Hopf link fibration. We will see more general examples in
section~\ref{S:Cvx2Ccv}.

\subsection{Between convexity and concavity}

As mentioned earlier, in the process of changing a convex boundary to
a concave boundary by attaching handles along a nicely fibered link we
encounter boundaries which are partially convex and partially
concave. Understanding how to control symplectic forms along
such boundaries is essential to the construction.
\begin{defn} \label{D:DCP} 
A {\em di\-la\-tion-con\-trac\-tion pair} on a symplectic $4$-mani\-fold
$(X,\omega)$ is a pair $(V^+,V^-)$ of vector fields defined
respectively on (possibly empty) open subsets $X^+$ and $X^-$ of $X$
such that the following equations hold (on the sets where they make
sense):
\[
  \mathcal{L}_{V^\pm} (\omega) = \pm \omega, \quad
  \omega(V^+,V^-) = 0\;
\]
A {\em contact pair} on a $3$-manifold $M$ is
a pair $(\alpha^+,\alpha^-)$ of 1-forms, defined respectively on
(possibly empty) open subsets $M^\pm$ of $M$, such that $M = M^+ \cup
M^-$ and such that the following equations hold (on the sets where
they make sense):
\[ \pm \alpha^\pm \wedge d \alpha^\pm > 0, \quad - d
  \alpha^- = d \alpha^+ 
\]
If $M$ is an oriented $3$-dimensional submanifold of $(X,\omega)$ then
$(V^+,V^-)$ {\em transversely covers} $M$ if $M \subset X^+ \cup X^-$
and if each vector field is positively transverse, where defined, to
$M$.  Letting $M^\pm = M \cap X^\pm$ and $\alpha^\pm =
\imath_{V^\pm}(\omega)|_{M^\pm}$, we see that the {\em induced pair}
$(\alpha^+,\alpha^-)$ is a contact pair on $M$.
\end{defn}
Notice that, for a contact pair $(\alpha^+,\alpha^-)$, $\alpha^+$ is a
positive contact form on $M^+$ while $\alpha^-$ is a negative contact
form on $M^-$. Together the two 1-forms give a globally defined,
closed, nondegenerate $2$-form $\gamma$ such that $\gamma|_{M^\pm} = \pm
d\alpha^\pm$. 

By a boundary which is partially convex and partially concave, we mean
a boundary which is transversely covered by a di\-la\-tion-con\-trac\-tion
pair. (Both convex boundaries and concave boundaries are special
cases.) By a germ of a symplectic form along a $3$-manifold $M$ we mean
an equivalence class of symplectic $4$-manifolds containing $M$, where
$(X_1,\omega_1) \sim (X_2,\omega_2)$ if there exist neighborhoods
$N_i$ of $M$ in $X_i$ and a symplectomorphism from $N_1$ to $N_2$
restricting to the identity on $M$.  We will prove the following
result in section~\ref{S:CvxCcvBdries}:
\begin{propn} \label{P:CPairSGerm}
A contact pair $(\alpha^+,\alpha^-)$ on a $3$-manifold $M$ defines a
unique symplectic germ $\mathcal{G}(\alpha^+,\alpha^-)$ along $M$ in
the following sense: 
\begin{enumerate} 
\item 
There exists a symplectic $4$-manifold $(X,\omega)$ containing $M$ with
a di\-la\-tion-con\-trac\-tion pair transversely covering $M$ inducing
$(\alpha^+,\alpha^-)$.
\item 
Any other symplectic $4$-manifold $(X_1,\omega_1)$ containing $M$ with
the property that $\omega_1 |_M = \pm d\alpha^\pm$ represents the
same germ along $M$.
\end{enumerate}
\end{propn}
This in particular implies that the induced contact pair on a
partially convex and partially concave boundary uniquely determines
the germ of the symplectic form along the boundary; we have already
seen this in the purely convex and purely concave cases.

\subsection{Outline}
We will now prove these results in the following order: After
establishing some terminology regarding handles, we will prove
theorem~\ref{T:WeakCvx} and discuss the relationship with the
Legendrian $2$-handles of Weinstein. (The constructions of the $2$-handles
in both theorems will be closely modelled on Weinstein's
construction.) Then we will prove the necessary results on partially
convex and partially concave boundaries and construct a class of
symplectic $2$-handles for this setting. Using this more general setup
we will prove theorem~\ref{T:cvx2ccv} and construct some
examples.

For background on basic tools in symplectic and contact constructions,
especially various versions of Darboux's theorem for symplectic and
contact structures, the reader is referred to~\cite{McDSal}.

\section{Terminology for Handles}
\label{S:ConvTerm}

\begin{figure}
\begin{center}
\includegraphics{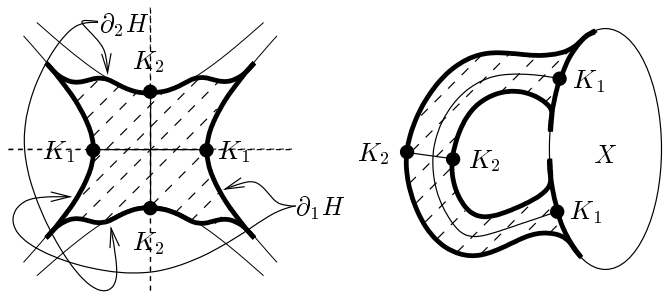}
\caption{Schematic pictures of a handle $H$ and of the result of attaching
$H$ to a manifold $X$}
\label{F:BasicHandle}
\end{center}\end{figure}

Our standard model for an $n$-dimensional $k$-handle will be a subset
$H$ of $\reals^n$ constructed in the following manner: Let $f$ be a
Morse function on $\reals^n$ with a single critical point of index $k$
at $0$, with $f(0)=0$. Choose constants $\epsilon_1 < 0 < \epsilon_2$;
$H$ will be a subset of $f^{-1}[\epsilon_1,\epsilon_2]$ bounded by two
smooth codimension one submanifolds with boundary, the ``attaching
boundary'' $\partial_1 H$ and the ``free boundary'' $\partial_2
H$. See figure~\ref{F:BasicHandle}. The attaching boundary $\partial_1
H$ is a closed tubular neighborhood of the descending sphere $K_1$ in
$f^{-1}\{\epsilon_1\}$, so that $\partial_1 H \cong S^{k-1} \times
B^{n-k}$. The free boundary $\partial_2 H$ begins as a tubular
neighborhood of the ascending sphere $K_2$ in $f^{-1}\{\epsilon_2\}$
but then dips down to join $\partial_1 H$ in $f^{-1}\{\epsilon_1\}$ so
that $\partial(\partial_1 H) = \partial(\partial_2 H)$. (Thus
$\partial_2 H \cong B^k \times S^{n-k-1}$.) Some form of smooth
``interpolation'' from $f^{-1}\{\epsilon_1\}$ to
$f^{-1}\{\epsilon_2\}$ must be specified to construct $\partial_2 H$;
in our constructions we will use a vector field transverse to the
level sets of $f$ to guide this interpolation.

Our convention is to orient $\partial_1 H$ and $\partial_2 H$ as level
sets of $f$, since $H$ is not technically a manifold with
boundary. This gives $\partial_2 H$ the orientation we would expect if
it really were a boundary component of $H$, but gives $\partial_1 H$
the nonstandard orientation; the two orientations agree on $\partial_1
H \cap \partial_2 H$. The descending sphere $K_1$ comes with a
canonical framing in $\partial_1 H$ which we will call the ``handle
framing'' of $K_1$.

Note in the figure that $\partial_1 H \cap \partial_2 H$ is a
codimension one submanifold, diffeomorphic to $S^{k-1} \times
S^{n-k-1} \times I$. This ``flange'' on the handle guarantees the
smoothness of the ``corners'' after attaching $H$ to an $n$-manifold
$X$. To attach $H$ to $X$, let $\hat{X}$ be $X$ with an open collar
attached to $\partial X$. We must specify an embedding of an open
neighborhood of $\partial_1 H$ in $\reals^n$ into $\hat{X}$
restricting to an embedding of $\partial_1 H$ into $\partial X$. Of
course, to do this we actually only need to specify the embedding of
$\partial_1 H$ into $\partial X$ but when we add symplectic structures
we should be more careful. In fact in the smooth case we need only
specify the image $K$ of $K_1$ in $\partial X$ and the framing $F$ of
$K$ induced by the handle framing of $K_1$ in $\partial_1 H$ to
completely determine the diffeomorphism type of the result of
attaching $H$ to $X$.

To construct a symplectic handle $(H,\omega_0)$ and attach it
symplectically to a symplectic manifold $(X,\omega)$ requires a
symplectic structure $\omega_0$ on $\reals^n$, an extension of
$\omega$ from $X$ to $\hat{X}$, and a symplectic embedding of a
neighborhood of $\partial_1 H$ into $\hat{X}$ restricting to an
embedding of $\partial_1 H$ into $\partial X$. If we have a symplectic
dilation positively transverse to $\partial_1 H$ and another
symplectic dilation positively transverse to $\partial X$, then, using
the symplectification of the contact forms, we need only specify an
embedding of $\partial_1 H$ into $\partial X$ which preserves the
induced contact forms in order to specify the symplectic embedding of
a neighborhood of $\partial_1 H$. If the embedding of $\partial_1 H$
into $\partial X$ instead only respects the contact structures, then
we must enlarge $(X,\omega)$ in a neighborhood of the image of
$\partial_1 H$ and perhaps rescale $\omega_0$ to arrange that the
embedding actually preserves the contact forms.

\section{Weakly convex boundaries and proof of theorem~\ref{T:WeakCvx}}
\label{S:WeakCvxBdries}

The idea in the proof of theorem~\ref{T:WeakCvx} is as follows: We use
polar coordinates on $\reals^4 = \reals^2 \times \reals^2$, by which
we mean coordinates $(r_1,\theta_1,r_2,\theta_2)$ where
$(r_i,\theta_i)$ are polar coordinates on the respective $\reals^2$
factors. We construct the $4$-dimensional $2$-handle $H$ as a neighborhood
of the origin in $\reals^4$, using the Morse function $f = -r_1^2 +
r_2^2$. We give $\reals^4$ the standard symplectic form and construct
a symplectic dilation $V$ which is transverse to the level sets of $f$
wherever it is defined, but which does not extend across
$\{r_1=0\}$. When we construct $H$, $V$ will be transverse to both
$\partial_1 H$ and $\partial_2 H \setminus K_2$, but will not extend
across $K_2$. Thus $V$ will induce positive contact structures on
$\partial_1 H$ and on $\partial_2 H \setminus K_2$. We will show that
the contact structure on $\partial_2 H$ can be deformed in a
neighborhood of $K_2$ so that it does extend across $K_2$, maintaining
the non-degeneracy condition needed to get weak convexity. Finally,
$K_1$ will be a transverse knot in $\partial_1 H$ and we will see why
a condition must be imposed on the framing in order to attach $H$
along a given transverse knot $K$.

\begin{proof}[Proof of theorem~\ref{T:WeakCvx}]
The standard symplectic form on $\reals^4$ is $\omega_0 = r_1 dr_1
d\theta_1 + r_2 dr_2 d\theta_2$.  Let $f = -r_1^2 + r_2^2$ and let $V
= \frac{1}{2} [ (r_1 - \frac{1}{r_1} ) \partial_{r_1} + r_2
\partial_{r_2} ]$. This vector field and some level sets of $f$
are shown in figure~\ref{F:weakVfield}.
\begin{figure} 
\begin{center}
\includegraphics{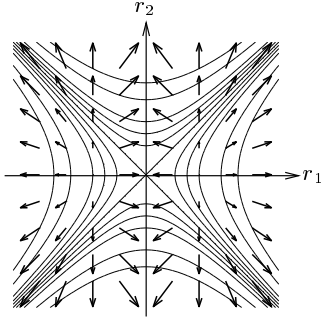}
\caption{Level sets of $f$ and the symplectic dilation $V$}
\label{F:weakVfield}
\end{center}\end{figure}

Notice that $V$ is a symplectic dilation and is positively transverse
to the level sets of $f$ as long as $f > -1$ and as long as $V$ is
defined, but that $V$ does not extend across $\{r_1 = 0\}$. Choose
constants $\epsilon_1$ and $\epsilon_2$ with $-1 < \epsilon_1 < 0 <
\epsilon_2$. The handle $H$ will be a subset of
$f^{-1}[\epsilon_1,\epsilon_2]$, with $\partial_1 H \subset
f^{-1}\{\epsilon_1\}$, $K_1 = \{r_2 = 0\} \cap
f^{-1}\{\epsilon_1\}$ and $K_2 = \{r_1 = 0\} \cap
f^{-1}\{\epsilon_2\}$.

First we  calculate the contact forms induced by $V$ on
$f^{-1}\{\epsilon_1\}$ and $f^{-1}\{\epsilon_2\} \setminus K_2$. Both
are restrictions of the form $\imath_V \omega_0 =
\half [(r_1^2 - 1) d\theta_1 + r_2^2 d\theta_2 ]$. On
$f^{-1}\{\epsilon_1\}$, natural polar coordinates consistent with our
orientation convention are
$(r=r_2,\mu=\theta_2, \lambda=-\theta_1)$.
With respect to these coordinates we get
$\alpha_1 =\half [r^2 d\mu - (r^2 - \epsilon_1 - 1) d\lambda]$.
Natural polar coordinates on $f^{-1}\{\epsilon_2\}$ are
$(r=r_1,\mu=\theta_1,\lambda=\theta_2)$ and
with respect to these coordinates we get
$\alpha_2 = \half [(r^2-1) d\mu + (r^2+\epsilon_2) d\lambda]$.

Now, given any small positive $\delta$, we will show how to modify the
contact structure $\xi_2 = \ker \alpha_2$ inside a neighborhood $\{r^2
\leq 1 +\delta\}$ of $K_2$ to get a new contact structure
$\xi_2^\prime$ which extends across $K_2$ and still satisfies the
property that $\omega_0|_{\xi_2^\prime}$ is nondegenerate. Let
$\xi_2^\prime = \ker \alpha_2^\prime$ where $\alpha_2^\prime =
d\lambda+t(r^2)d\mu$ and $t(r^2)$ goes smoothly to $0$ as $r^2$ goes
to $0$, has positive derivative and is equal to
$\frac{r^2-1}{r^2+\epsilon_2}$ on $\{r^2 \geq 1+\delta\}$. This agrees
with $\xi_2$ outside $\{r^2 \leq 1+\delta\}$ because $\xi_2 = \ker
(d\lambda+\frac{r^2-1}{r^2+\epsilon_2} d\mu)$, is contact because
$t^\prime > 0$ and satisfies the nondegeneracy condition because
$\alpha_2^\prime \wedge (\omega_0|_{f^{-1}\{\epsilon_2\}}) = (1-t(r^2))r
dr \wedge d\mu \wedge d\lambda > 0$. Intuitively, $\xi_2$ twists too
far as we move in towards $K_2$ to extend across $K_2$, so we back off
to $\{r^2 = 1+\delta\}$ and then twist more slowly so that
$\xi_2^\prime$ does extend.

Forward flow along $V$ for time $t$ starting at a point
$(r,\mu,\lambda) \in f^{-1}\{\epsilon_1\}$ gives a map $\Phi$ from a
subset of $\reals \times f^{-1}\{\epsilon_1\}$ into $\reals^4$ defined
by the following equations:
\[
\begin{array}{cc}
    r_1^2 \circ \Phi = (r^2 - \epsilon_1 - 1) e^t + 1\;, & 
    \theta_1 \circ \Phi = -\lambda \\
    r_2^2 \circ \Phi = r^2 e^t\;, &
    \theta_2 \circ \Phi = \mu \;
\end{array}
\]
Letting $R = \sqrt{\epsilon_2 (\frac{1+\epsilon_1}{1+\epsilon_2})}$
and $T=\log(\frac{1 + \epsilon_2}{1 + \epsilon_1})$, we see that
forward flow for time $t=T$ defines a diffeomorphism $\phi :
f^{-1}(\epsilon_1) \setminus \{ r \leq R \} \to f^{-1}(\epsilon_2)
\setminus K_2$. Note that $1+\epsilon_1 > R^2$ and that $\phi\{r^2 =
1+\epsilon_1\} = \{r^2=1\}$.

Given any three radii $R_3 > R_2 > R_1 > \sqrt{1+\epsilon_1}$ we can
construct a symplectic handle $H$ as follows: Choose a smooth function
$h:[0,R_3] \to [0,T]$ which is equal to $T$ on $[0,R_1]$, is
decreasing on $[R_1,R_2]$ and is equal to $0$ on $[R_2,R_3]$. Then let
$H$ be the union of all the forward flow lines starting at points $p
\in \{r \leq R_3\} \subset f^{-1}\{\epsilon_1\}$ flowing for time less
than or equal to $h(r(p))$, together with $\{r_1 = 0\} \cap
f^{-1}[\epsilon_1,\epsilon_2]$. The attaching boundary $\partial_1 H$
is $\{r \leq R_3\} \subset f^{-1}\{\epsilon_1\}$ while the free
boundary $\partial_2 H$ is the image under $\Phi$ of the graph of $h$
in $\reals \times f^{-1}\{\epsilon_1\}$ together with $K_2$. This
construction is illustrated in figure~\ref{F:WeakHandle}. In the
figure $R_3$, $R_2$ and $R_1$ are rather far apart for the sake of
clarity, but in general one would carry out this construction with
these radii only slightly larger than $\sqrt{1+\epsilon_1}$. A few
forward flow lines for $V$ are shown, starting on
$f^{-1}\{\epsilon_1\}$.

\begin{figure} 
\begin{center}
\includegraphics{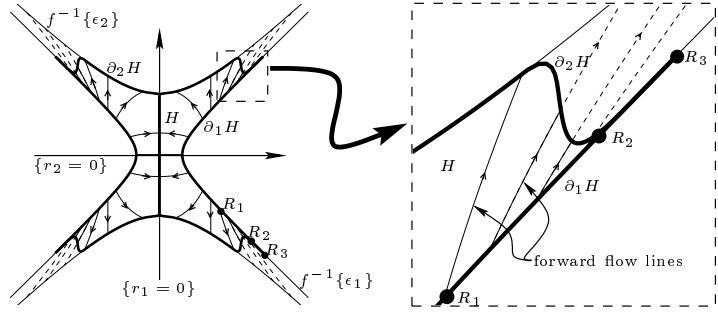}
\caption{Construction of $H$ in the proof of
theorem~\ref{T:WeakCvx}}
\label{F:WeakHandle}
\end{center}\end{figure}

The construction of $\xi_2^\prime$ on
$f^{-1}\{\epsilon_2\}$ gives a positive contact structure
$\xi_2^{\prime\prime}$ on $\partial_2 H$ by letting
$\xi_2^{\prime\prime} = \xi_2^\prime$ on $\partial_2 H \cap
f^{-1}\{\epsilon_2\}$ and letting $\xi_2^{\prime\prime} = \ker(
\imath_V \omega_0|_{\partial_2 H})$ elsewhere.  The nondegeneracy
condition on $\xi_2^\prime$ and the fact that $V$ is positively
transverse to $\partial_2 H$ where defined gives us that
$\omega_0|_{\xi_2^{\prime\prime}}$ is nondegenerate. Furthermore
$\xi_2^{\prime\prime}$ agrees with $\xi_1 = \ker \alpha_1$ on $\partial_1 H
\cap \partial_2 H$, so that if we succeed in attaching $H$ to
$(X,\omega)$ via a contactomorphic embedding of $(\partial_1 H, \xi_1)$ into
$(\partial X,\xi)$, the contact structures will patch together to give a
contact structure $\xi_Y$ on the boundary $\partial Y$ of the new 
symplectic manifold $(Y,\eta)$ with $\eta|_{\xi_Y}$ nondegenerate.

Now we show how to find this embedding and attach $H$.  Given the
transverse link $L \subset \partial X$ which is fat with respect to a
framing $F$, consider a component $K$ of $L$ with neighborhood $\nu$
with almost normal coordinates $(r,\mu,\lambda)$ which go out as far
as $F$. Notice that we can choose these coordinates so that $F_\mu =
F$. This means that the contact structure $\xi$ near $K$ is spanned by
$\partial_r$ and $\log(r^2) \partial_\mu + \partial_\lambda$ and that
$\{r \leq 1\} \subset \nu$. Choose some $\epsilon>0$ such that $\{r
\leq 1+\epsilon\} \subset \nu$. Now consider the coordinates
$(r,\mu,\lambda)$ on $\partial_1 H$; these are almost normal
coordinates and we could explicitly reparametrize $r$ to make them
normal. However it is sufficient to notice that $\xi_1$ is spanned by
$\partial_r$ and $\frac{r^2-\epsilon_1-1}{r^2} \partial_\mu +
\partial_\lambda$ and that this implies that, after reparametrizing
$r$, the coordinates would go out as far as the handle framing
$F_\mu$. This is because we construct $H$ with
$R_3>R_2>R_1>\sqrt{1+\epsilon_1}$. This means that if we construct $H$
with $R_3-\sqrt{1+\epsilon_1}$ small enough we can guarantee a
contactomorphism from $\partial_1 H$ to $\{r \leq 1+\epsilon\} \subset
\nu$, taking the handle framing to $F$.
\end{proof}

An attractive feature of this construction is that it is
straightforward to visualize the resulting ``contact surgery''. The
new contact manifold $(\partial Y, \xi_Y)$ is constructed from
$(\partial X, \xi)$ as follows: For each component $K_i$ of the
transverse link $L$, let $\nu_i$ be a neighborhood of $K_i$ with
normal coordinates $(r_i,\mu_i,\lambda_i)$ which go out at least as
far as the framing $F_i$. For each $i$ choose a small positive
$\epsilon_i$ so that $\{r_i \leq e^{n_i} + \epsilon_i\}$ is in the
image of the coordinate system, where $F_i = F_{\mu_i} + n_i$. For
each $i$ remove the solid torus $\{r_i < e^{n_i} + \epsilon_i/2\}$ and
glue back in a solid torus $\nu_i^\prime$ by overlapping along
$\{e^{n_i} + \epsilon_i/2 \leq r_i \leq e^{n_i}+\epsilon_i\}$,
realizing the (topological) surgery with framing $F_i$. Then we can
extend $\xi$ across $\nu_i^\prime$ exactly because, on $\{e^{n_i} +
\epsilon_i/2 \leq r_i \leq e^{n_i}+\epsilon_i\}$, $\xi$ is twisting
towards the longitude realizing the framing $F_i$, which, after the
surgery, will become a meridian. We needed to remove at least $\{r_i
\leq e^{n_i}\}$ because, when $r_i < e^{n_i}$, $\xi$ has already
twisted past this longitude.

To see that this construction achieves all the surgeries achievable
using Weinstein's Legendrian surgeries, we now present the proof of
proposition~\ref{P:LegXvrs}, beginning with a lemma.
\begin{lemma} \label{L:LegXvrs}
If $\beta_1$ and $\beta_2$ are 1-forms on disks $D_1$ and $D_2$ such
that $d\beta_i > 0$ then $\alpha_i = d\lambda + \beta_i$ are positive
contact forms on $D_i \times S^1$ ($\lambda$ being the
$S^1$-coordinate).  If $\phi: (D_1,d\beta_1) \to (D_2,d\beta_2)$ is a
symplectomorphism then there exists a function $h: D_1 \to \reals$
such that the diffeomorphism $\Phi:D_1 \times S^1 \to D_2 \times S^1$
given by $\Phi(p,\lambda) = (\phi(p),\lambda + h(p))$ satisfies
$\Phi^* \alpha_2 = \alpha_1$.
\end{lemma}
\begin{proof}
That each $\alpha_i$ is a positive contact form is a straightforward
calculation. To see why the rest is true, note that since $\phi^*
d\beta_2 = d\beta_1$, the 1-form $\beta_1 - \phi^* \beta_2$ is closed
and therefore exact. Choose $h$ so that $\beta_1 = \phi^* \beta_2 +
dh$. Then $\Phi^* \alpha_2 = d(\Phi^* \lambda) + \phi^*\beta_2 =
d\lambda + dh + \phi^* \beta_2 = \alpha_1$.
\end{proof}

\begin{proof}[Proof of proposition~\ref{P:LegXvrs}]
Let $K$ be the Legendrian knot, $\nu$ a neighborhood of $K$ and $F
\leq \tb(K)-1$ a framing of $K$.  Without loss of generality, by
Darboux's theorem for contact structures, we may assume that $\nu$ has
the form $ (\nu = D_\epsilon \times S^1, \xi = \ker \alpha)$ where
$\alpha = dy - x d\lambda$, $x$ and $y$ are coordinates on the disk
$D_\epsilon = \{x^2 + y^2 < \epsilon^2\}$ of radius $\epsilon$, and
$\lambda$ is the $S^1$-coordinate. This is because $K = \{0\} \times
S^1$ is Legendrian. Note that in this model $\tb(K)$ is the
``zero-framing'' coming from the product structure on $\nu$. We will
measure framings relative to this product framing, so that $\tb(K)-1 =
-1$.

On $\{ x > 0 \}$, we have $\xi = \ker \alpha^\prime$, where
$\alpha^\prime = d\lambda + \beta$ and $\beta = -\frac{1}{x} dy$. The
$2$-form $d\beta = \frac{1}{x^2} dx \wedge dy$ is positive on
$\{x>0\}$. We will construct a symplectomorphism $\phi$ from the disk
$D_2$ of radius $2$ in $\reals^2$ with the standard symplectic form
$dx \wedge dy = r dr \wedge d\mu$ onto a region $D \subset \{ x > 0,
x^2 + y^2 < \epsilon \}$ with the symplectic form $d\beta$. Then on
$D_2$ let $\beta_2 = \half r^2 d\mu$ and note that $d\beta_2 = r dr
\wedge d\mu$. Thus lemma~\ref{L:LegXvrs} gives a contactomorphism
$\Phi$ from $(D_2 \times S^1, d\lambda+\beta_2)$ onto $(D \times S^1,
\alpha^\prime)$ taking the zero framing to the zero framing.  Finally
note that $K_2 = \{r=0\} \subset (D_2 \times S^1, d\lambda + \beta_2)$
is fat with respect to the framing $-1 = \tb(K)-1$.

We construct $\phi$ directly. Choose two positive constants $c_1$
and $c_2$, define $\phi$ by:
\[ \phi(x,y) = (\frac{c_2}{c_1 -  x}, c_2 y) \]
and verify that $\phi^* \frac{1}{x^2} dx\wedge dy = dx \wedge dy$. The
map is only defined when $x < c_1$, but as long as we choose $c_1 >
2$, $\phi$ will be defined on $D_2$. By choosing $c_2$ small enough we
can guarantee that $\phi(D_2) \subset \{x >0,x^2 + y^2 < \epsilon\} =
D$.
\end{proof}  

\section{Partially convex and partially concave boundaries}
\label{S:CvxCcvBdries} 

Before proving theorem~\ref{T:cvx2ccv} we need to develop a theory of
symplectic boundaries which are partially convex and partially
concave. (We only develop this theory in dimension four.) The
definitions were given in section~\ref{S:ResMot}; our first task is to
show that a contact pair induced by a di\-la\-tion-con\-trac\-tion pair
uniquely determines the germ of the symplectic form and to look at
some corollaries of the proof. After that we will show how to
construct a symplectic $2$-handle $H$ with a di\-la\-tion-con\-trac\-tion pair
transversely covering $\partial_1 H$ and $\partial_2 H$ and inducing
specified contact pairs. In the next section we will use these tools
to prove theorem~\ref{T:cvx2ccv}.

\subsection{Uniqueness of germs determined by contact pairs}

Henceforth the notation $(M,(\alpha^+,\alpha^-))$ will refer to a
$3$-manifold equipped with a contact pair. We will always refer to the
domains of the forms as $M^\pm$ and will let $M^0 = M^+ \cap M^-$ with
$\alpha^0 = \alpha^+|_{M^0} + \alpha^-|_{M^0}$. A positive
(resp. negative) contact form is a special case of a contact pair,
with $M^-=\emptyset$ (resp. $M^+ = \emptyset$).  Let $R_{\alpha^\pm}$
be the Reeb vector fields for $\alpha^\pm$, and note that $\alpha^0$
is closed and nowhere zero and that $\alpha^0(R_{\alpha^\pm}) > 1$.

\begin{lemma} \label{L:ExUnDCPs}
Given $(M,(\alpha^+,\alpha^-))$,
consider the two symplectic manifolds $(S^+, \omega^+)$ and $(S^-,
\omega^-)$, where $S^\pm = \reals \times M^\pm$ and $\omega^\pm =
\pm d (e^{\pm t} \alpha^\pm)$, and identify $M^\pm$ with $\{0\}
\times M^\pm \subset S^\pm$.  Firstly, if $(X,\omega)$ is
another symplectic manifold containing $M^+$ or $M^-$ with a single
symplectic dilation $V^+$ or contraction $V^-$ positively transverse
to $M^\pm$ inducing the contact form $\alpha^+$ or $\alpha^-$,
respectively, then flow along $V^\pm$ starting from $M^\pm$ gives an
embedding $\Phi$ of an open subset of $S^+$ or $S^-$, respectively,
into $X$ such that $\Phi^* \omega = \omega^\pm$ and
$D\Phi(\partial_t) = V^\pm$.

Secondly: 
\begin{enumerate}
\item 
There exists a unique vector field $V^-$ defined on $\reals \times
M^0$ such that the pair $(V^+ = \partial_t, V^-)$ is a
di\-la\-tion-con\-trac\-tion pair on $(S^+,\omega^+)$ inducing the contact
pair $(\alpha^+, \alpha^-|_{M^0})$ on $M^+$.
\item 
There exists a unique vector field $V^+$ defined on $\reals \times
M^0$ such that the pair $(V^+, V^- = \partial_t)$ is a
di\-la\-tion-con\-trac\-tion pair on $(S^-,\omega^-)$ inducing the contact
pair $(\alpha^+|_{M^0}, \alpha^-)$ on $M^-$.
\end{enumerate}
\end{lemma}

We will call the symplectic manifold $(S^+,\omega^+)$ with its
di\-la\-tion-con\-trac\-tion pair the {\em positive symplectification} of
$(M^+,(\alpha^+,\alpha^-|_{M^0}))$ and we will call $(S^-,\omega^-)$
with its di\-la\-tion-con\-trac\-tion pair the {\em negative
symplectification} of $(M^-,(\alpha^+|_{M^0},\alpha^-))$.

\begin{proof}
The first result follows from the fact that $\Phi^* \omega$ and
$\omega^\pm$ are solutions to the same ordinary differential equations
with the same initial conditions.

For the second result, let $\gamma = \pm d\alpha^\pm$, let $g^\pm =
\frac{\alpha^0 \wedge \gamma}{\alpha^\pm \wedge \gamma} =
\alpha^0(R_{\alpha^\pm})$ and let $\beta^\pm = \alpha^0 - g^\pm
\alpha^\pm$. We will show that there exist unique vector fields $Z^\pm
\in \ker \alpha^+ \cap \ker \alpha^-$ on $M^0$ such that
$\imath_{Z^\pm}(\gamma) = \beta^\pm$. Then the vector fields $V^-$ and
$V^+$, on $(S^+,\omega^+)$ and $(S^-,\omega^-)$ respectively, are
given by the following formulae:
\begin{enumerate}
\item $V^- = (g^+ e^{-t} -1) \partial_t + e^{-t} Z^+$
\item $V^+ = (g^- e^t - 1) \partial_t + e^t Z^-$
\end{enumerate}
We first show the existence and uniqueness of $Z^\pm$, then show that
$V^\pm$ as described in these formulae satisfy the conditions of the
lemma, and then show uniqueness of $V^\pm$.
  
There exists a unique $Z^\pm \in \ker \alpha^\pm$ such that
$\imath_{Z^\pm}(\gamma) = \beta^\pm$ because contraction with $\gamma$
gives a linear isomorphism from $\ker \alpha^\pm$ to $\{\beta \mid
\beta \wedge \gamma = 0 \}$ (this depends on working in dimension
$3$), and $\beta^\pm$ is constructed to be in this latter
subspace. But $Z^\pm$ is also in $\ker \alpha^0$ because $0 =
\gamma(Z^\pm,Z^\pm) = \beta^\pm(Z^\pm) = \alpha^0(Z^\pm)$, and thus
$Z^\pm \in \ker \alpha^\mp$.

On $S^+$, letting $V^+ = \partial_t$ and $V^- = (g^+ e^{-t} - 1)
\partial_t + e^{-t} Z^+$, we need to show that
\begin{align}
\mathcal{L}_{V^\pm} (\omega^+) &= \pm \omega^+ \label{E:LieDer} \\
\imath_{V^\pm}(\omega^+) |_{t=0} &= \alpha^\pm \label{E:Contr} \\
\omega^+(V^+,V^-) &= 0 \label{E:Lagr}
\end{align}
First note that $\omega^+ = e^t(dt \wedge \alpha^+ + \gamma)$.
Equation~\ref{E:Lagr} is quick: $\omega^+(V^+,V^-) = e^{-t}
\omega^+(\partial_t,Z^+) = \alpha^+(Z^+) = 0$. To show
equation~\ref{E:LieDer} and equation~\ref{E:Contr}, note that
$\imath_{\partial_t}(\omega^+) = e^t \alpha^+$ and that
\begin{align*}
\imath_{V^-}(\omega^+) &= (g^+ e^{-t} - 1) \imath_{\partial_t}
(\omega^+) + e^{-t} \imath_{Z^+} (\omega^+) \\
&= e^t[(g^+ e^{-t} - 1) \alpha^+ + e^{-t} \beta^+] \\
&= - e^t \alpha^+ + \alpha^0
\end{align*}
  
Next we will prove that $V^-$ is the unique vector field on $\reals
\times M^0$ satisfying these equations. Suppose $V^-_0$ and $V^-_1$
are two solutions. Let $\delta = \imath_{V^-_1 - V^-_0} (\omega^+)$;
we will prove that $\delta = 0$ and thus conclude that $V^-_0 =
V^-_1$. Equation~\ref{E:LieDer} implies that $d \delta = 0$,
equation~\ref{E:Contr} implies that $\delta|_{\{t=0\}} = 0$ and
equation~\ref{E:Lagr} implies that $\delta(\partial_t) = 0$
everywhere. Thus $\delta$ is invariant in the $t$ direction and
vanishes when $t=0$, so $\delta = 0$ everywhere.

On $S^-$ the argument is a mirror image of the argument for $S^+$.
\end{proof}

Note that the uniqueness argument first proved uniqueness {\em along
$M$} (since $\delta|_{\{t=0\}} = 0$) and then proved
uniqueness for all $t$. Thus in fact we have also proved:
\begin{lemma} \label{L:UnDCPAlongM}
In $(\reals \times M^0, \omega^+)$ there exists a unique vector field
$V^-$ along $M^0$ (i.e. a section of $T_{M^0} (\reals \times M^0)$)
positively transverse to $M^0$ such that $\imath_{V^-} \omega^+ =
\alpha^-$ and $\omega^+(\partial_t,V^-) = 0$. Likewise there exists a
unique vector field $V^+$ along $M^0$ in $(\reals \times M^0,
\omega^-)$ positively transverse to $M^0$ such that $\imath_{V^+}
\omega^- = \alpha^+$ and $\omega^-(V^+,\partial_t) = 0$.
\end{lemma}

We will make use of the following observation in later constructions:
\begin{lemma} \label{L:BothTransverse}
Given $(M,(\alpha^+,\alpha^-))$, let $(S^\pm,\omega^\pm)$ be the
positive and negative symplectifications with the
di\-la\-tion-con\-trac\-tion pairs from lemma~\ref{L:ExUnDCPs}.  Let
$g^\pm$, $\beta^\pm$ and $Z^\pm$ be as in the proof of
lemma~\ref{L:ExUnDCPs}. Given a function $h: M^\pm \to \reals$,
consider its graph $M^\pm_h = \{(h(p),p)\} \subset S^\pm$. Then $V^\pm
= \partial_t$ is automatically positively transverse to $M^\pm_h$ and
$V^\mp$ is positively transverse to $M^\pm_h$ if and only if $e^{\pm
h} < g^\pm - dh(Z^\pm)$ on $M^0$. Identifying $M_h$ with $M$ in the
obvious way, the induced contact pair on $M_h$ is then given by
\[
\alpha^\pm_h = e^h \alpha^\pm,\; \alpha^0_h = \alpha^0,\;
\alpha^\mp_h = \alpha^0 - e^h \alpha^\pm\; .
\]
\end{lemma}
\begin{proof}
Everything follows from the explicit expressions $V^\mp = (g^\pm
e^{\mp t} -1) \partial_t + e^{\mp t} Z^\pm$ and $\omega^\pm = e^{\pm
t}(dt \wedge \alpha^\pm + \gamma)$.
\end{proof}

\begin{proof}[Proof of Proposition~\ref{P:CPairSGerm}]
First we prove existence of $(X,\omega)$ with its
di\-la\-tion\--con\-trac\-tion pair $(V^+,V^-)$.  Construct
$(S^\pm,\omega^\pm)$ as in lemma~\ref{L:ExUnDCPs}. Let
$U^+$ be an open neighborhood of $M^0$ in $S^+$ such that flow along
$V^+$ in $S^-$ starting from $p_0 \in M^0$ is defined for all times
$t$ with $(t,p_0) \in U^+$.  This gives an embedding $\phi^+: U^+
\into S^-$ such that $\phi^+|_{M^0} = \id$ and $D\phi^+(\partial_t) =
V^+$. Since both $\partial_t$ and $V^+$ induce the same contact forms
on $M^0$ and are both symplectic dilations, we can conclude that
$(\phi^+)^*(\omega^-) = \omega^+$. By the uniqueness in
lemma~\ref{L:ExUnDCPs} we also know that $D\phi^+(V^-) = \partial_t$,
so that $\phi^+$ preserves all the relevant structure. Let $U^- =
\phi^+(U^+)$ and $\phi^- = (\phi^+)^{-1}$.
  
Now choose two functions $f^\pm: M \to [0,\infty)$ such that
$f^\pm|_{M^\mp \setminus M^0} = 0$ but $f^+ + f^- > 0$
everywhere, let $F^\pm = \{(t,p) \mid -f^\pm(p) < t < f^\pm(p)
\} \subset S^\pm$ and let $E^\pm = F^\pm \cap \phi^\mp(F^\mp
\cap U^\mp)$.  Finally let $\psi^\pm = \phi^\pm|_{E^\pm} :
E^\pm \to E^\mp$.  If we choose $f^\pm$ small enough we can
guarantee that
\[ 
X = F^+ \cup_{\psi^+} F^- 
\]
is Hausdorff. Since $(\psi^+)^* \omega^- = \omega^+$ and
$D\psi^+(V^\pm) = V^\pm$, we know that the symplectic forms and
the di\-la\-tion-con\-trac\-tion pairs patch together to define a
symplectic form $\omega$ on $X$ with a di\-la\-tion-con\-trac\-tion
pair $(V^+,V^-)$ transversely covering $M$ inducing
$(\alpha^+,\alpha^-)$. (See figure~\ref{F:Sympation}.)
\begin{figure} 
\begin{center}
\includegraphics{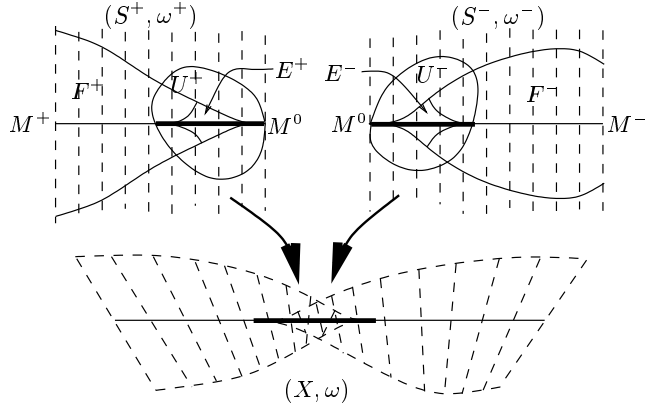}
\caption{Glueing $S^+$ to $S^-$ after some trimming} 
\label{F:Sympation}
\end{center}
\end{figure}
  
For the uniqueness result, we need to construct a bundle isomorphism
$\Psi : T_M X_1 \to T_M X$ covering the identity and preserving the
symplectic forms, and then we can apply Darboux's theorem. To do this
we construct a pair of vector fields $(V^+_1,V^-_1)$ along $M$ in $T_M
X_1$ with open domains $M^\pm_1 \subset M^\pm$ covering $M$, both
positively transverse to $M$, such that $\imath_{V^\pm_1}
\omega_1|_{M^\pm_1} = \alpha^\pm|_{M^\pm_1}$ and such that
$\omega_1(V^+_1,V^-_1) = 0$.  Then there exists a unique $\Psi$
sending $V^\pm_1$ to $V^\pm|_{M_1^\pm}$ by lemma~\ref{L:UnDCPAlongM}.
To see the existence of the pair $(V^+_1,V^-_1)$, first extend
$\alpha^+$ to a maximal rank 1-form on $T_{M^+} X_1$ to get $V^+_0$
along $M^+$ such that $\imath_{V^+_0} \omega_1 |_{M^+} =
\alpha^+$. Then by lemma~\ref{L:UnDCPAlongM} there exists a unique
$V^-_0$ along $M^0$ such that $\imath_{V^-_0} \omega_1 |_{M^0} =
\alpha^-|_{M^0}$ and $\omega_1(V^+_0,V^-_0) = 0$. Let
$\tilde{\alpha}^-_0 = \imath_{V^-_0} \omega_1$, a maximal rank 1-form
on $T_{M^0} X_1$ which extends $\alpha^-|_{M^0}$. Then there exists a
maximal rank extension $\tilde{\alpha}^-_1$ of $\alpha^-$ to $T_{M^-}
X_1$ which agrees with $\tilde{\alpha}^-_0$ outside a closed set $C$
inside $M^-$ containing $M^- \setminus M^0$ in its interior. Let
$V^-_1$ be the corresponding vector field such that $\imath_{V^-_1}
\omega_1 = \tilde{\alpha}^-_1$ and let $V^+_1 = V^+_0 |_{M^+ \setminus
C}$.

\end{proof}

\begin{cor} \label{C:DiffPairsSameGerm}
Suppose that $(X,\omega)$ is a symplectic $4$-manifold with $\partial X$
partially convex and partially concave, with induced contact pair
$(\alpha_1^+,\alpha_1^-)$. Notice that $\omega|_{\partial X} = \pm
d\alpha_1^\pm$.  Now suppose that $(\alpha_2^+,\alpha_2^-)$
is another contact pair on $\partial X$ with $\omega|_{\partial X} = \pm
d\alpha_2^\pm$. Then $\mathcal{G}(\alpha_1^+,\alpha_1^-) =
\mathcal{G}(\alpha_2^+,\alpha_2^-)$, so that $\partial X$ can also be
seen as partially convex and partially concave with induced contact
pair $(\alpha_2^+,\alpha_2^-)$. This holds in particular if
$(\alpha_2^+,\alpha_2^-)$ is obtained from $(\alpha_1^+,\alpha_1^-)$
by extending or restricting the domains of the 1-forms maintaining the
contact pair properties.
\end{cor}

\subsection{Construction of a handle}

Now we will show how to construct a symplectic $2$-handle $H$
according to the following plan: As in section~\ref{S:WeakCvxBdries},
$H$ will be a subset of $f^{-1}[\epsilon_1,\epsilon_2] \subset
(\reals^4,\omega_0)$ for the Morse function $f=-r_1^2 + r_2^2$, where
$\epsilon_1 < 0 < \epsilon_2$. In this case, however, we will consider
a particular di\-la\-tion-con\-trac\-tion pair $(V^+,V^-)$ on
$\reals^4$ with $V^+$ defined on $\reals^4 \setminus \{r_1=0\}$ and
$V^-$ on $\reals^4 \setminus \{r_2=0\}$. This pair will transversely
cover both $\partial_1 H$ and $\partial_2 H$, inducing contact pairs
$(\alpha^+_1,\alpha^-_1)$ and $(\alpha^+_2,\alpha^-_2)$. Let
$\partial_i^\pm H$ be the respective domains of these 1-forms; we then
have $\partial_1^+ H = \partial_1 H$, $\partial_1^- H = \partial_1 H
\setminus K_1$, $\partial_2^+ H = \partial_2 H \setminus K_2$ and
$\partial_2^- H = \partial_2 H$. We will use the flow along $V^+$ to
guide the construction of $\partial_2 H$. The tricky point is to
arrange that {\em both} vector fields end up transverse to $\partial_2
H$. For this we will use lemma~\ref{L:BothTransverse} and choose our
``interpolation'' function $h$ carefully.

We use the same coordinates on $f^{-1}\{\epsilon_1\}$ and
$f^{-1}\{\epsilon_2\}$ that we used in section~\ref{S:WeakCvxBdries}:
On $f^{-1}\{\epsilon_1\}$, these are $(r=r_2,\mu=\theta_2,
\lambda=-\theta_1)$ while on $f^{-1}\{\epsilon_2\}$ these are
$(r=r_1,\mu=\theta_1,\lambda=\theta_2)$.

To characterize the contact pairs on $\partial_1 H$ we need a little
more terminology. Let $\nu$ be a neighborhood of a knot $K$ with a
contact pair $(\alpha^+,\alpha^-)$ and with polar coordinates
$(r,\mu,\lambda)$. Let the domain of $\alpha^\pm$ be $\nu^\pm$, assume
that $\nu^+ = \nu$ and $\nu^- = \nu \setminus K$ and that
$(r,\mu,\lambda)$ is an almost normal coordinate system for $\ker
\alpha^+$. Let $\alpha^0 = \alpha^+ + \alpha^-$ (in which case the
pair $(\alpha^+,\alpha^0)$ determines $\alpha^-$).
\begin{defn} \label{D:PrepSurg}
Such a contact pair $(\alpha^+,\alpha^-)$ is {\em well-behaved} with
respect to $(r,\mu,\lambda)$ if $R_{\alpha^+} = A
\partial_\mu + B \partial_\lambda$ and $\alpha^0 = C d \mu + D
d\lambda$ for constants $A,B,C,D$ with $B$ and $C$ positive. The
contact pair is {\em prepared for surgery} with respect to
$(r,\mu,\lambda)$ if we also have that $A=B$, $D>0$ and $AD > 1$. In
each case, we may also say that the knot $K$ is well-behaved or
prepared for surgery.
\end{defn}
The fact that $\alpha^0(R_{\alpha^+}) > 1$ implies that $AC+BD>1$. We
will call the quadruple $(A,B,C,D)$ the {\em structural data} for the
contact pair (with respect to the coordinate system); this data
completely determines $(\alpha^+,\alpha^-)$ on $\nu$ up to a
reparametrization of $r$. A convenient model, given $(A,B,C,D)$, is
$\alpha^+ = \frac{1}{B+Ar^2}(r^2 d\mu + d\lambda)$ and $\alpha^- =
\alpha^0 - \alpha^+ = C d\mu + D d\lambda - \alpha^+$. (Simply verify
that this is a well-behaved contact pair with the given data.)
\begin{propn} \label{P:CPHandle}
Given a contact pair $(\alpha^+,\alpha^-)$ on a neighborhood
$\nu$ of a knot $K$ which is prepared for surgery with structural data
$(A=B,C,D)$, there exists a handle $H$ such that the induced pair
$(\alpha_1^+,\alpha_1^-)$ on $\partial_1 H$ is prepared for surgery
with the same structural data (with respect to the above coordinates
on $\partial_1 H \subset f^{-1}\{\epsilon_1\}$). Furthermore,
with the contact pair on $f^{-1}\{\epsilon_1\}$ fixed, we can
construct $H$ so as to make $\partial_1 H$  arbitrarily small as a
neighborhood of $K_1$.
\end{propn}
Thus, given any symplectic $4$-manifold $X$ with boundary which is
partially convex and partially concave with induced contact pair
$(\alpha^+,\alpha^-)$, we can attach such a handle to $X$
whenever we can find a knot $K$ in $(\partial X, (\alpha^+,\alpha^-))$
which is prepared for surgery with respect to some coordinate
system. The coordinate system determines the framing, and the handle
can be attached along an arbitrarily small neighborhood of $K$. This
is in contrast to the construction in section~\ref{S:WeakCvxBdries},
which required the attaching neighborhood to be ``fat enough'' with
respect to the chosen framing.
\begin{proof}
Let $\omega_0 = r_1 dr_1 d\theta_1 + r_2 dr_2 d\theta_2$ on
$\reals^4$ with $f=-r_1^2 + r_2^2$, and let
\begin{align*}
  V^+ &= (\half r_1 - \frac{C}{r_1}) \partial_{r_1} + \half r_2
  \partial_{r_2} \\ V^- &= -\half r_1 \partial_{r_1} - (\half r_2 -
  \frac{D}{r_2}) \partial_{r_2} 
\end{align*}

Calculation shows that $(V^+,V^-)$ is a di\-la\-tion-con\-trac\-tion pair
which transversely covers the level sets of $f$ as long as $-2D < f <
2C$. Let $\epsilon_1 = \frac{2}{A} - 2D$ and note that $-2D <
\epsilon_1 < 0$ (because $A>0$ and $AD > 1$). Choose any $\epsilon_2$
with $0 < \epsilon_2 < 2C$.  Then the induced contact pair on
$f^{-1}(\epsilon_1)$ is given by:
\[ 
\alpha^+_1 = \half [r^2 (d\mu -
d\lambda)] + \frac{1}{A} d\lambda, \quad \alpha^0_1 = C d\mu + D
d\lambda, \]
with
\[ \quad R_{\alpha^+_1} = A (\partial_\mu +
\partial_\lambda)\;.
\] 
Forward flow along $V^+$ gives a map $\Phi$ from some subset of
$\reals \times f^{-1}(\epsilon_1)$ into $\reals^4$ defined by:
\[
\begin{array}{cc}
    r_1^2 \circ \Phi = (r^2 -\frac{2}{A})e^t + 2D\;, & 
    \theta_1 \circ \Phi = -\lambda \\
    r_2^2 \circ \Phi = r^2 e^t\;, &
    \theta_2 \circ \Phi = \mu \;
\end{array}
\]
Letting $R=\frac{\epsilon_2}{A(D+\epsilon_2/2)}$ and $T = \log( A(D +
\frac{\epsilon_2}{2}))$ (which is positive because $AD > 1$), we see
that forward flow for time $t = T$ defines a diffeomorphism
$\phi:f^{-1}\{\epsilon_1\} \setminus \{ r^2 \leq R\} \to
f^{-1}\{\epsilon_2\} \setminus K_2$.  We can make $R$ arbitrarily
small by choosing $\epsilon_2$ small enough. Now for any choice of
radii $R_3 > R_2 > R_1 > R$ we can build $H$ exactly as in
section~\ref{S:WeakCvxBdries}, using a function $h:[0,R_3] \to [0,T]$
which is $T$ on $[0,R_1]$, decreasing on $[R_1,R_2]$ and $0$ on
$[R_2,R_3]$.

By construction $\partial_2 H$ is transverse to $V^+$. The only part
of $\partial_2 H$ which could fail to be transverse to $V^-$ is $
\Gamma_h = \Phi(\{(h(r(p)),p) \mid R_1 \leq r(p) \leq R_2 \})$.  Using
lemma~\ref{L:BothTransverse} we can state conditions for $\Gamma_h$ to
be transverse to $V^-$. Using the notation from the proof of
lemma~\ref{L:ExUnDCPs} we get:
\begin{align*}
\gamma &= r dr \wedge (d\mu - d\lambda) \\
g^+ &= A(C+D) \\
\beta^+ &= (C-A(C+D)r^2) (d\mu - d\lambda) \\
Z^+ &= \frac{1}{r}(C-A(C+D)r^2) \partial_r 
\end{align*}
Thus $h$ needs to satisfy
\[ 
e^{h(r)} < g^+ - h^\prime(r) dr(Z^+) =
A(C+D) - \frac{1}{r} h^\prime(r) (C-A(C+D)r^2) \;
\] 
for all $r \in [R_1,R_2]$.  Since we will have $h^\prime \leq 0$, $e^h
\leq A(D+\frac{\epsilon_2}{2})$ and $\epsilon_2 < 2C$, this will work
as long as
\[
r^2 < \frac{C}{A(C+D)}\; .
\]
Therefore if we build $H$ with $R_2 < \sqrt{\frac{C}{A(C+D)}}$ we can
guarantee that $\partial_2 H$ is transverse to both $V^+$ and
$V^-$. We see that $R_1$, $R_2$ and $R_3$ can be chosen arbitrarily
small, as long as we choose $\epsilon_2$ small enough, so that we can
arrange for $\partial_1 H$ to be an arbitrarily small neighborhood of
$K_1$.
\end{proof}

\subsection{Preparing well-behaved knots for surgery}
In order to prove theorem~\ref{T:cvx2ccv} we will want to attach
handles along well-behaved transverse knots, so now we present a
method to turn well-behaved knots into knots which are prepared for
surgery under certain conditions

Recall that, given a neighborhood $\nu$ of a knot $K$ with contact
structure $\xi$, if $(r,\mu,\lambda)$ is an almost normal coordinate
system on $\nu$ then so is $(r,\mu - k \lambda, \lambda)$ for any $k
\in \mathbb{Z}$. If $\xi = \ker \alpha^+$ and $(\alpha^+,\alpha^-)$ is
well-behaved with respect to $(r,\mu,\lambda)$ then
$(\alpha^+,\alpha^-)$ will also be well-behaved with respect to
$(r,\mu-k\lambda,\lambda)$, and the structural data $(A,B,C,D)$ with
respect to $(r,\mu,\lambda)$ will change to $(A-Bk,B,C,D+Ck)$ with
respect to $(r,\mu-k\lambda,\lambda)$. Thus, if we are willing to
increase framings, it is easy to arrange that $D > 0$ and that $BD >
1$. However it is not clear how to arrange that $A=B$.
\begin{lemma} \label{L:WellB2PrepSurg}
Suppose that $(\alpha^+,\alpha^-)$ is well-behaved with respect to
coordinates $(r,\mu,\lambda)$ on $\nu$ with structural data
$(A,B,C,D)$, where $D>0$ and $BD > 1$.  Then, for any $\epsilon > 0$,
there exists a function $h:\nu \to [0,\infty)$ supported inside $\{r
\leq \epsilon\}$ with the following properties: Let $S^+=(\reals \time
\nu,\omega^+)$ be the positive symplectification of
$(\nu,(\alpha^+,\alpha^-))$, with di\-la\-tion-con\-trac\-tion pair
$(V^+,V^-)$. Let $\nu_h = \{(h(p),p)\} \subset S^+$ and let $\pi:\nu_h
\to \nu$ be the natural projection. Then $\nu_h$ is transverse to both
vector fields, and the induced contact pair $(\alpha^+_h,\alpha^-_h)$
is prepared for surgery with respect to the coordinates
$(r\circ\pi,\mu \circ \pi, \lambda \circ \pi)$ on some (smaller)
neighborhood of $\pi^{-1}(K)$.
\end{lemma}

\begin{proof}
As mentioned earlier we may assume that $\alpha^+ = \frac{1}{B + Ar^2}
(r^2 d\mu + d\lambda)$. To avoid too much notation, we will use
$(r,\mu,\lambda)$ on $\nu_h$ to refer to $(r\circ \pi, \mu \circ \pi,
\lambda \circ \pi)$. For a given $h$ we will have $\alpha^+_h = e^{h}
\alpha^+$ on $\nu_h$ and lemma~\ref{L:BothTransverse} gives us conditions for
$V^-$ to be transverse to $\nu_h$. When $V^-$ is also transverse to
$\nu_h$ we get $\alpha^0_h = \alpha^0$. Choose a constant $A_0$ with
$\frac{1}{D} < A_0 < B$ (we have $\frac{1}{D} < B$ because $BD>1$ and
$D>0$). Note that we then have $A_0 (C+D)>1$
(because $C>0$ and $D>0$). We will construct $h$ so that $\alpha^+_h =
\frac{1}{A_0 + A_0 r^2} (r^2 d\mu + d\lambda)$ on $\{r \leq \delta\}
\subset \nu_h$ for some positive $\delta < \epsilon$; this together
with $\alpha^0_h = \alpha^0$ gives a contact pair on $\nu_h$ which is
prepared for surgery on $\{r \leq \delta\}$.

Notice that this in fact determines $h$ on $[0,\delta]$ because we
must have $e^{h} = \frac{B+Ar^2}{A_0(1+r^2)}$ for $r \in
[0,\delta]$. We should check that $h$ so defined is in fact
positive: $h \geq 0$ as long as $B+Ar^2 \geq A_0(1+r^2)$, which
will hold for small enough $r$ as long as $B > A_0$, which is how we
chose $A_0$. In other words, if we choose $\delta$ small enough we can
guarantee that $h > 0$ on $\{ r \leq \delta\}$.

Next we check that, with $h$ thus defined on $\{r \leq \delta\}$,
$V^-$ is transverse to $\nu_h$. Calculating and applying
lemma~\ref{L:BothTransverse} we see that transversality will hold if 
\begin{equation}
e^{h} < AC + BD - \frac{(C - Dr^2)(B + Ar^2)}{2r} \frac{\partial
h}{\partial r} \label{E:xvrs}
\end{equation}
which, for our given $h$, becomes:
\[
B+Ar^2 < A_0(D+C)(B+Ar^2)\;.
\]
This holds for $r \leq \delta$ because $B+Ar^2 \geq A_0(1+r^2)$ and
$A_0(D+C) > 1$.

Now we should check that we can extend $h$ to $\nu$ to have support
inside $\{r \leq \epsilon^\prime\}$ for some $\epsilon^\prime <
\epsilon$, in such a way that $V^-$ remains transverse to $\nu_h$. On
$\{r \geq \epsilon^\prime\}$ the transversality condition~\ref{E:xvrs}
above will be satisfied because $h$ will be identically $0$ and $AC+BD
> 1$. For $r \leq \epsilon^\prime$, if we choose $\epsilon^\prime$
small enough we can replace condition~\ref{E:xvrs} by the following
simpler condition:
\begin{equation}
e^h < AC+BD - CB \frac{\partial h}{\partial r}
\label{E:xvrs2}
\end{equation}
Using the facts that $AC +BD >1$ and $C$ and $B$ are positive,
it is easy to extend $h$ to $r \leq \epsilon^\prime$ maintaining this
condition if $\frac{\partial h}{\partial r}(\delta) \leq 0$. If
$\frac{\partial h}{\partial r}(\delta) > 0$ then, after perhaps making
$\delta$ smaller still, we extend $h$ to $\{r \leq
\epsilon^\prime\}$, with $h=0$ near $\epsilon^\prime$, in such a way
that $\frac{\partial h}{\partial r}(r) < \frac{\partial h}{\partial
r}(\delta)$ for all $r > \delta$ and that $\frac{\partial h}{\partial
r}(r) <0$ for all $r > \delta + \delta_1$, for some small $\delta_1 >
0$. This is enough to conclude that conditon~\ref{E:xvrs2} is met fir
$r \leq \epsilon^\prime$.
\end{proof}

\begin{cor} \label{C:WellBEnlarge}
Suppose that $(X,\omega)$ is a symplectic $4$-manifold with $\partial X$
partially convex and partially concave with induced contact pair
$(\alpha^+,\alpha^-)$, that $K \subset \partial X$ is a transverse
knot with a neighborhood $\nu$ with coordinates $(r,\mu,\lambda)$ and
that $(\alpha^+,\alpha^-)$ is well-behaved with respect to
$(r,\mu,\lambda)$ on $\nu$ with structural data $(A,B,C,D)$. If $D>0$
and $BD>1$ then we can enlarge $(X,\omega)$ inside $\nu$ so as to
arrange that $\partial X$ is partially convex and partially concave
with induced contact pair $(\alpha^+,\alpha^-)$ which is prepared for
surgery on $\nu$ with respect to $(r,\mu,\lambda)$. Then we can attach
a handle as in proposition~\ref{P:CPHandle} along $K$ with framing
$F_\mu$.
\end{cor}
\begin{proof}
Enlarge $(X,\omega)$ using the positive symplectification of
$(\alpha^+|_\nu,\alpha^-|_\nu)$ (see 
lem\-ma~\ref{L:ExUnDCPs}).  Attach the subset $\{(t,p) \mid 0 \leq t
\leq h(p)\}$ of this positive symplectification, where $h$ comes from
lemma~\ref{L:WellB2PrepSurg}, using the uniqueness of the symplectic
germ $\mathcal{G}(\alpha^+|_\nu,\alpha^-|_\nu)$.
\end{proof}

\section{From convexity to concavity via fibered links}
\label{S:Cvx2Ccv}

In theorem~\ref{T:cvx2ccv}, after attaching the handles, $\partial Y$
comes with a link $L^\prime$, the union of the ascending circles of
the handles.  We will now prove theorem~\ref{T:cvx2ccv} and along the
way see the following characterization of the induced negative contact
form $\alpha_Y^-$ on $\partial Y$.
\begin{add}[to Theorem~\ref{T:cvx2ccv}] \label{A:cvx2ccv}
There exists a closed tubular neighborhood $\tau$ of $L$, a
constant $k$ and a diffeomorphism $\phi$ from $\partial X \setminus
\tau$ to $\partial_Y \setminus L^\prime$ such that $k dp -
\phi^*(\alpha^-_Y) = e^{h} \alpha$ for some function $h: \partial X
\setminus \tau \to [0,\infty)$.
\end{add}
\begin{proof}
[Proof of Theorem~\ref{T:cvx2ccv} and addendum~\ref{A:cvx2ccv}] We
will first argue that we can enlarge $(X,\omega)$ so as to arrange
that $\partial X$ is in fact partially convex and partially concave,
with induced contact pair $(\alpha^+,\alpha^-)$, and so as to arrange
that there are coordinates near each component of $L$, realizing the
desired framing, with respect to which $(\alpha^+,\alpha^-)$ is
well-behaved satisfying the conditions of
corollary~\ref{C:WellBEnlarge}.

Recall the notation in the definition of ``nicely fibered''. The
transverse contact vector field is $V$ and the fibration is $p:
\partial X \setminus L \to S^1$. With $(X,\omega)$ as given, we have
some induced contact form $\alpha$ on $\partial X$ with $\xi = \ker
\alpha$ such that $(\partial X,\xi,L,p)$ satisfy the definition of
``nicely fibered''. Consider a new contact form $\alpha^+$
defined by $\alpha^+|_\xi = 0$ and $\alpha^+(V) = 1$. Notice
that we then have $R_{\alpha^+} = V$.  The new contact form
$\alpha^+$ has the same kernel as $\alpha$, so $\alpha^+ = g
\alpha$ for some function $g:\partial X \to \reals \setminus
\{0\}$. By replacing $V$ with $-V$ if necessary we can arrange that $g
> 0$. Also note that we can replace $V$ with $kV$ for any constant
$k>0$ without changing the ``nicely fibered'' condition, and thus we
can arrange that $g > 1$ (using the compactness of $\partial X$). This
means that we can enlarge $(X,\omega)$ so as to arrange that the
induced contact form on $\partial X$ is actually $\alpha^+$,
using the symplectification $(\reals \times \partial X, d(e^t
\alpha))$. 

Now choose a constant $c$ so that $c \cdot dp(V) > 1$ on $M \setminus
L$ (this depends on the compactness of $\partial X$ and on the fact
that $V$ and $dp$ are invariant on a neighborhood of $L$ so that
$dp(V)$ is constant near $L$). Let $\alpha^0 = c \cdot dp$ and let
$\alpha^- = \alpha^0 - \alpha^+$. Then $(\alpha^+,\alpha^-)$ is a
contact pair on $\partial X$ with $-d\alpha^- = d\alpha^+$, so using
corollary~\ref{C:DiffPairsSameGerm} we may now regard $\partial X$ as
partially concave and partially convex with induced contact pair
$(\alpha^+,\alpha^-)$. Of course $\alpha^+$ is still defined on all of
$\partial X$ while $\alpha^-$ is only defined on $\partial X \setminus
L$, and so for now $\alpha^-$ contains no new information. However,
when we show that $(\alpha^+,\alpha^-)$ is well-behaved near $L$ we
will be able to attach the handles from the previous section, after
which $\alpha^-$ will extend across the new boundary and we will be
able to forget about $\alpha^+$ to conclude that the new boundary is
concave. Thus $\alpha^-$ contains the seed of the concavity which we
will achieve after attaching handles.

We see that $(\alpha^+,\alpha^-)$ is well-behaved near $L$ using the
``nicely fibered'' condition again. Near each component we know that
there are normal coordinates $(r,\mu,\lambda)$ such that
$\alpha^0(\partial_r) = 0$, $dr(V) = 0$ and such that $V$ and
$\alpha^0$ are invariant under the flows of $\partial_r$,
$\partial_\mu$ and $\partial_\lambda$.  This immediately establishes
that $\alpha^0 = C d\mu + D d\lambda$ for some constants $C$ and $D$,
and that $R_{\alpha^+} = A \partial_\mu + B \partial_\lambda$ for some
constants $A$ and $B$. Now we need to arrange that $B$ and $C$ are
positive; this will follow from the orientation condition on the
characteristic foliation on the fibers.  Looking at our orientation
convention, the fact that the foliation points radially inwards means
that $\alpha^+ \wedge \alpha^0 \wedge (-dr) > 0$. Recall that we can
take $\alpha^+ = \frac{1}{B+Ar^2} (r^2 d\mu+d\lambda)$ as our model,
so we get that $\frac{C-Dr^2}{B+Ar^2} > 0$ for small $r$. Thus either
$B$ and $C$ are both positive, or if they are both negative we can
replace the coordinate system with $(r,-\mu,-\lambda)$ to get $B$ and
$C$ both positive.

Now if we arrange that each coordinate system $(r,\mu,\lambda)$ for
each component of $L$ realizes the desired framing $F$ of $L$, we see
that the condition that $F$ is positive with respect to the fibration
$p$ means exactly that, near each component, $-\frac{D}{C} < 0$ and
since $C>0$ this means that $D>0$. Now we still may not have that
$BD>1$. To arrange this we may need to again replace $\alpha^0$ by $k
\alpha^0$ for some constant $k>1$ (again using compactness of
$\partial X$), which will replace $C$ and $D$ with $kC$ and $kD$. Now
corollary~\ref{C:WellBEnlarge} shows how to enlarge $(X,\omega)$ and
attach handles.  After enlarging $(X,\omega)$, we have $\partial X$
partially convex and partially concave with induced contact pair
$(\alpha^+,\alpha^-)$, with $\partial^+ X = \partial X$ and
$\partial^- X = \partial X \setminus L$. After attaching the handles
$\partial Y$ is partially convex and partially concave with induced
contact pair $(\alpha^+_Y, \alpha^-_Y)$ with domains $\partial^\pm Y$,
with $\partial^+ Y = \partial Y \setminus L^\prime$ and $\partial^- Y
= \partial Y$ (where $L^\prime$ is the union of the ascending
spheres). This means that we can ignore $\alpha^+_Y$ and realize that
in fact $\partial Y$ is concave with induced negative contact form
$\alpha^-_Y$, and the characterization of $\alpha^-_Y$ in the addendum
follows.
\end{proof}

\subsection{Examples}

First we will show that $S^3$ with surgery on either the unknot
with any framing $F>0$ or the Hopf link with any framing $F \geq
0$ can be realized as the concave boundary of a symplectic
$4$-manifold. 

Let $(r_1,\theta_1,r_2,\theta_2)$ be polar coordinates on $\reals^4$.
Consider $S^3 = \partial B^4 \subset (\reals^4, \omega)$ where $\omega
= r_1 dr_1 \wedge d\theta_1 + r_2 dr_2 \wedge d\theta_2$. Then $V =
\half(r_1 \partial_{r_1} + r_2 \partial_{r_2})$ is a symplectic
dilation transverse to $S^3$ inducing the standard positive contact
form $\alpha = \half(r_1^2 d\theta_1 + r_2^2 d\theta_2)$ on $S^3$. We
compute that $R_\alpha = \partial_{\theta_1} + \partial_{\theta_2}$
and let $V=R_\alpha$, a transverse contact vector field. Now consider
two cases:
\begin{enumerate}
\item 
$L = K = \{r_1 = 0\} \subset S^3$, the standard unknot.  Consider the
fibration $p = \theta_1: S^3 \setminus L \to S^1$ and notice that
$dp(V) >0$. Polar coordinates near $K$ are given by
$(r=r_1,\mu=\theta_1,\lambda=\theta_2)$, from which we can verify that
$K$ is nicely fibered with fibration $p$ and contact vector field
$V$. With respect to these coordinates, $F_\mu$ is the standard
$0$-framing of $K$, and since $dp = d\mu$, the condition that a
framing $F$ is positive with respective to $p$ is simply the condition
that $F > 0$. Thus we can attach a single ``convex-to-concave'' handle
along the unknot with any framing $F > 0$ to make a symplectic
manifold with concave boundary.
\item 
$L = K_1 \cup K_2$ where $K_i = \{ r_i = 0 \} \subset S^3$, the
standard Hopf Link. Now consider the fibration $p = \theta_1 +
\theta_2 : S^3 \setminus L \to S^1$ and again, $dp(V) >0$. Oriented
polar coordinates near $K_1$ are as above and oriented polar
coordinates near $K_2$ are given by
$(r=r_2,\mu=\theta_2,\lambda=\theta_1)$ and again we verify that $L$
is nicely fibered with fibration $p$ and contact vector field
$V$. Also $F_\mu$ is the standard $0$-framing of each $K_i$. Now,
however, $dp = d\mu + d\lambda$ near each $K_i$, so the condition that
a framing $F$ is positive with respect to $p$ is the condition that
the framing on each $K_i$ is greater than $-1$.  Thus we can attach a
pair of ``convex-to-concave'' handles along the Hopf link as long as
each handle is framed with framing $0$ or larger, and the result is a
symplectic manifold with concave boundary.
\end{enumerate}

The first example generalizes. Given an $n$-punctured surface $\Sigma$
with a proper Morse function $f_\Sigma: \Sigma \to [0,\infty)$ with
only critical points of index $0$ and $1$, with the critical points of
index $0$ lying in $f^{-1}\{0\}$ and those of index $1$ lying in
$f^{-1}\{\frac{1}{4}\}$, we can use Weinstein's construction in
dimension $2$ to get a symplectic form $\omega_\Sigma$ on $\Sigma$ and
a gradient-like symplectic dilation $V_\Sigma$ such that the structure
is ``standard'' on, say, $f^{-1}(\half,\infty)$. In other words,
$f^{-1}(\half,\infty)$ looks like $n$ copies of $\reals^2 \setminus
\{r^2 \leq \half\}$ with its standard symplectic form $r dr \wedge
d\theta$, symplectic dilation $\half r \partial_r$ and Morse function
$r^2$.  Consider the symplectic $4$-manifold $(\Sigma \times \reals^2,
\omega = \omega_\Sigma + r dr \wedge d\theta)$ with the Morse function
$f = f_\Sigma + r^2$. Note that $V = V_\Sigma + \half r \partial_r$ is
a gradient-like symplectic dilation for $f$. Let $X = f^{-1}[0,1]$ and
let $M = \partial X = f^{-1}(1)$; $M$ is the convex boundary of
$(X,\omega)$ with induced contact form $\alpha = \alpha_\Sigma + \half
r^2 d\theta$, where $\alpha_\Sigma = \imath_{V_\Sigma}
\omega_\Sigma$. Note that $M$ is diffeomorphic to the ``boundary with
smoothed corners'' of the product of a disk and a compact surface of
genus $g$ with $n$ boundary components. Alternately, if the handle
decomposition of $\Sigma$ is the usual one with one 0-handle and
$(2g+n-1)$ 1-handles, then $X$ is diffeomorphic to $B^4$ with
$(2g+n-1)$ 1-handles attached, and $M$ is diffeomorphic to $S^3$ with
$0$-surgery on $(2g+n-1)$ unlinked unknots.

We can decompose $M$ into two open sets: $A = \{ f_\Sigma < 1, r^2 = 1
- f_\Sigma \}$ and $B = \{0 \leq r^2 < \half, \half < f_\Sigma = 1 -
r^2 \leq 1 \}$. The set $A$ is the complement of the $n$-component
link $L = \{ r^2 = 0, f_\Sigma = 1 \}$, and the function $\theta: A =
M \setminus L \to S^1$ is a fibration with fiber diffeomorphic to
$\Sigma$. To see that $d\theta (R_\alpha) > 0$, note that this is
equivalent to the requirement that $d\theta \wedge d\alpha > 0$. But
$d\theta \wedge d\alpha = d\theta \wedge d\alpha_\Sigma = d\theta
\wedge \omega_\Sigma > 0$.  On $B$, the entire structure is identical
to $n$ copies of the structure described in the earlier example on
$S^3$ in a neighborhood of the standard unknot. Thus we can attach
handles along $L$ with any framing larger than the ``$0$-framing''
determined by the fibration to create a concave symplectic $4$-manifold.

To build on these constructions it is worth noting a simple variation
of Weinstein's construction in~\cite{Weinstein}: Suppose $(X,\omega)$
is a symplectic $4$-manifold with {\em concave} boundary and $K \subset
\partial X$ is a Legendrian knot with respect to the induced negative
contact structure on $\partial X$. Then a symplectic handle can be
attached along $K$ with framing $\tb(K)-1$ such that the new
symplectic manifold again has {\em concave} boundary. This can be seen
as follows: Weinstein's ``convex-to-convex'' $4$-dimensional $2$-handles
are constructed as subsets of $\reals^4$ using the symplectic form
$\omega = dx_1 \wedge dy_1 + dx_2 \wedge dy_2$, the Morse function
$f=-x_1^2- x_2^2 + y_1^2 + y_2^2$ and the symplectic dilation $V= -x_1
\partial_{x_1} + 2 y_1 \partial_{y_1} - x_2 \partial_{x_2} + 2 y_2
\partial_{y_2}$. One checks that $V$ is positively transverse to level
sets of $f$ and that the descending sphere $K_1 \subset
f^{-1}\{\epsilon_1\}$ (for any $\epsilon_1 < 0$) is Legendrian with
respect to the induced contact structure on $f^{-1}\{\epsilon_1\}$ and
that the handle framing of $K_1$ is one less than the
Thurston-Bennequin framing, so that such a handle can be attached
along any Legendrian knot $K$ with framing $\tb(K)-1$. To construct
``concave-to-concave'' $2$-handles, use the same symplectic form and
Morse function but now consider the symplectic {\em contraction} $V =
- 2 x_1 \partial_{x_1} + y_1 \partial_{y_1} - 2 x_2 \partial_{x_2} +
y_2 \partial_{y_2}$. Now simply check that this new $V$ is positively
transverse to level sets of $f$, that $K_1$ is Legendrian with respect
to the induced negative contact structure on $f^{-1}\{\epsilon_1\}$
and that the handle framing is again equal to $\tb(K_1)-1$.

\begin{figure}
\begin{center}
\includegraphics{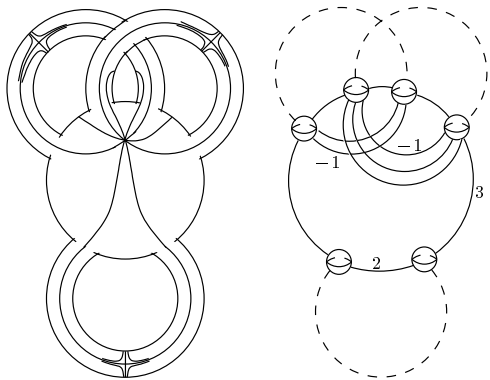}
\caption{A handlebody decomposition of a surface with some leaves of
the characteristic foliation giving a framed link description of a
symplectic $4$-manifold with concave boundary.}
\label{F:samplelink}
\end{center}
\end{figure}

In our example above ($X = f^{-1}[0,1] \subset \Sigma \times \reals^2$
and $(M,\alpha) = (\partial X, \alpha_\Sigma + \half r^2 d\theta)$),
the characteristic foliation on each fiber $\Sigma$ is given by the
flow lines of the vector field $V_\Sigma$.  We can construct $\Sigma$
so that there will be some closed leaves of this singular foliation
(containing singular points). Each such closed leaf is a Legendrian
knot in $(M,\alpha)$ the Thurston-Bennequin framing of which is the
framing given by the fiber. Let $Y$ be the result of attaching
$2$-handles along $L=\{r^2=0,f_\Sigma=1\}$, with $\partial Y$ concave
with induced negative contact form $\alpha_Y$ and fibered link
$L^\prime$. Notice that, under the diffeomorphism $\phi: \partial X
\setminus \tau \to \partial Y \setminus L^\prime$ of
addendum~\ref{A:cvx2ccv}, $\alpha_Y$ induces the same characteristic
foliation on the fibers of $\theta$ as $\alpha$ did. Thus the closed
leaves of the foliation are again Legendrian knots in $(\partial Y,
\alpha_Y)$ and we can attach symplectic $2$-handles along these knots
(as in the previous paragraph) with framing $-1$ with respect to the
fibers to build larger manifolds with concave boundary. This increases
our class of examples of $3$-manifolds which bound concave symplectic
$4$-manifolds; first perform any positive surgeries along the original
$n$-component link $L$, then perform $-1$ surgeries on arbitrarily
many copies of each closed leaf of the singular foliation (use
different fibers to get the different copies). A random example of
this construction is shown in figure~\ref{F:samplelink}.

These examples are fairly immediate and we hope that more
sophisticated examples can be constructed using these techniques. The
obvious challenge is to construct examples in which the negative
contact structure on the concave boundary is recognizable as
contactomorphic via an orientation reversing diffeomorphism to the
positive contact structure on some other convex boundary, so that
interesting closed symplectic manifolds can be constructed.

\providecommand{\bysame}{\leavevmode\hbox to3em{\hrulefill}\thinspace}

\end{document}